\documentstyle[pb-diagram]{article}

\newcommand{\ba}{\begin{array}}
\newcommand{\ea}{\end{array}}

\newtheorem{theorem}{THEOREM}[section]

\newtheorem{corollary}{COROLLARY}[section]
\newtheorem{definition}{DEFINITION}[section]
\newtheorem{proposition}{PROPOSITION}[section]

\newcommand{\be}{\begin{enumerate}}
\newcommand{\ee}{\end{enumerate}}
\newcommand{\bi}{\begin{itemize}}
\newcommand{\ei}{\end{itemize}}
\newcommand{\bd}{\begin{description}}
\newcommand{\ed}{\end{description}}

\newcommand{\et}{\wedge}

\newcommand{\vel}{\vee}

\newcommand{\imp}{\rightarrow}

\newcommand{\beq}{\begin{eqnarray*}}
\newcommand{\eeq}{\end{eqnarray*}}

\author{ {\bf F. Parlamento,  }
\\Department of Mathematics and Computer Science
\\University of Udine,  via  Delle Scienze 206, 33100 Udine, Italy.\\
e-mail: {\em  franco.parlamento$@$uniud.it}}
 \title{Truth-value semantics and functional extensions for  classical  logic of partial terms based on equality
 \footnote{Work supported by funds PRIN-/MIUR. The author is grateful to the referee for very helpful comments and suggestions.} }
\date{}

\begin{document}

\maketitle

\begin{abstract}

We develop a  bottom-up approach to truth-value semantics for classical logic of partial terms based on equality, and apply it to prove the conservativity  of the addition of partial description and selection functions, independently of any strictness assumption.

\end{abstract}


\subsection{Introduction}

We assume the reader is familiar with the natural deduction system  
for classical first order logic,
conceived   as the result of the direct analysis of actual mathematical reasoning, as presented by Gentzen in \cite{G35}. 
At the same time we ask her/him to leave aside, for a moment, the now standard classical set theoretic formulation of the notion of logical consequence.
By classical logic of partial terms based on equality we mean the standard natural deduction system,
with the proviso, of a semantical nature,  that not all terms are assumed to be necessarily denoting;  a  feature that is syntactically reflected by the restriction of
 the  usual $\forall$-elimination and $\exists$-introduction rules, as formulated in \cite{P65},
  to  variables or individual parameters only.
On the other hand,  that  a term $t$ is  denoting  is expressed by the assumption $\exists x (x=t)$, for $x$ not occurring in $t$, in agreement with Quine's Thesis\footnote{So christened in \cite{H59} and expressed by Quine's  dictum  from \cite{Q48}, ``to be is to be the value of a variable".}, as originally proposed in 
  \cite{H59} and \cite{LH59}.
 Truth-valued semantics has  been extensively investigated by H. Leblanc, among others, see {\cite{L71}, \cite{L76} and especially \cite{L01},
which presents it as     the result of a progressive simplification of the standard set theoretic semantics, first to countable models, then to  Henkin's models and finally to no model at all.
Quite to the opposite, we wish to show that  truth-value semantics can be approached from below, so to speak, by following the search of the simplest mathematical means by which one can establish that a proposition is not deducible from others, by the application of the given natural deduction  rules, if that is indeed the case. 
We  will explain to what extent that approach   determines the usual truth tables for the propositional connectives and how it  leads to truth-value semantics, when quantifiers are involved.
A distinguished feature of our treatment, with respect  to Leblanc's,  is that it deals with first order languages endowed with function symbols, which, apart from its intrinsic interest, is clearly necessary if $tv$-semantics for partial logic has to be applied to show the conservativity of the addition of partial description and selection functions.
As in  \cite{G36} and \cite{P65},  we refer
to  the articulation of a first order language in which, beyond a countable supply of variables,   meant to be used for quantification, one has  also an infinite supply of individual parameters, meant to remain free names for generic objects of whatever (non empty) domain one happens to be talking about.
 Once truth-value  semantics   ($tv$-semantics, for short) is defined, we will sketch a proof that it  is indeed fully adequate, namely that not only our motivating goal, namely correctness,
  but also completeness holds.
Then  we  establish  the {\em Extension Property}, which will be basic for all later developments. The basic idea to deal semantically  with the undefinedness of a pure term $t$ with respect to a truth-value valuation ($tv$-valuation for short) $v$, is simply to say that $t$ is non denoting with respect to $v$ if for all individual parameters $a$, $v(a=t)=\bf f$.
  Our main purpose is then to employ $tv$-semantics  to show that the above logical framework  is  appropriate
  to deal with non empty domains, with a language in which individual parameters stand for objects of the domain but more general terms, such as $-1$ or $1/(a-a)$, when the natural  or the real numbers are involved, need not denote any object whatsoever.
See \cite{F95} for a more extended and very illuminating discussion. In fact, by using $tv$-semantics, we will prove   the conservativity of the  addition of partial selection and description functions, also when   the strictness axioms,  to the effect that: $1)$ all  constants are denoting, $2)$ if $ft_1\ldots t_n$ is denoting, then $t_1, \ldots, t_n$ are  denoting as well and $3)$ for $p$ other than $=$, if $p t_1\ldots t_n$ holds, then $t_1, \ldots, t_n$ are denoting, are added to the underlying logical framework.
 To obtain our conservativity results, we have obviously to take into account all possible $tv$-valuations: those for which there is a non denoting term can be disposed with by choosing one such term. For the remaining  ones, to be called {\em  totally denoting $tv$-valuations}, we have to enrich the language with a new constant: the {\em undefined} $\uparrow$, and show that the given valuation can be extended to the new language in a way that actually leaves $\uparrow$ undefined.
 To deal with the strictness axioms,  we have to adopt a  corresponding type of $tv$-valuation  and show that the Extension Property applies to them as well.
 The conservativity of the addition of partial selection functions  and partial description functions, with or without strictness axioms, then follows by a straightforward correctness/completeness  argument.
 Finally
  it is to be noted that totally denoting valuations are elementarily equivalent to   classical set theoretic structures (with total functions interpreting function symbols) and strict valuations are elementarily equivalent to  set theoretic structures with partial functions  interpreting function symbols.
  As  such,  totally denoting $tv$-valuations constitute  a natural intermediate step for the introduction  of what has become the standard semantics for classical first order logic, with completeness achieved as  a simple corollary. Correctness, on the other hand, crucially depends on proving the substitution lemmas
 (which, presumably, involves the tedious details mentioned in Gumb's obituary of Leblanc \cite{G00}). \footnote{ 
 Leblanc found truth-value semantics to be a useful teaching device enabling students to grasp more easily fundamental semantic concepts, because it abstracted from tedious details in standard, set-theoretic semantics. }

 \subsection{Pure terms and formulae}
 
 \begin{definition} Given a first order language ${\cal L}$,
 \bi
 \item[a)] A term $t$ of ${\cal L}$ is pure if  no variable occurs in $t$,
  \item[b)] A formula  $F$ of ${\cal L}$ is pure if  no variable occurs free in $F$,
 \ei
 The collection of pure terms of ${\cal L}$ will be denoted by $PureTerm_{\cal L}$.
  \end{definition}
  
 In particular, sentences are pure formulae. This terminology is inspired by  Gentzen's suggestion in \cite{G35} \footnote{``rein logische Formel"} and in \cite{S69} p. 70.\footnote{
 The concept of a formula is ordinarily used in a more general sense; the special case defined [above] might thus perhaps described as a Òpurely logical formulaÓ.}

\subsection{Natural deduction systems for partial logic}

As for the deductive apparatus we refer to the natural deduction system, which we denote 
by $N_c$, in which
the $\forall$-elimination and $\exists$-introduction rule take the restricted  form

\[
\ba{ccc}
\forall x F&~~~& F\{ x / {\bf y} \} \\ 
\cline{1-1} \cline{3-3}
F\{ x/ {\bf y}\} &   &\exists x F
\ea
\]
where ${\bf y}$ is either a free variable or an individual parameter.
A deduction is said to be pure when it involves pure formulae only,
in particular in its $\forall$-elimination and $\exists$-introduction,  ${\bf y}$ must be a parameter.
$G_1, \ldots, G_n \rhd_c F$ denotes that there is a  deduction in $N_c$ with conclusion $F$ and active assumptions included among $G_1,\ldots, G_n$.

\subsection{ A "bottom-up" approach to truth-value semantics}

At the propositional level, when required to explain why, for example,  $A$ does not follow from $A\imp B$ and $B$, one usually provides examples taken from the ordinary
or mathematical language, like   letting $A$ be ``the car runs out of gas" and $B$ be ``the car stops", where all is relevant is  our persuasion that if $A$ is true then $B$ is true as well, but if $B$ is true $A$ need not  necessarily be true. That naturally leads to the idea of a valuation of the propositional  atoms of the propositions we are investigating,  into at least two values.
Our goal of showing  that $F$ does not follow from $G_1, \ldots, G_n$ is reached if:

\bi
\item a method of computing values for compound statements is found such that one specific value, say ${\bf t}$,  is preserved by deductions, and a valuation $v$ of the propositional atoms in $G_1, \ldots, G_n, F$ is found, such that $G_1, \ldots, G_n$ takes the value ${\bf t}$, but $F$ does not.
\ei
Clearly for that to work at least two values are needed. Classical propositional semantics makes the minimal choice of two values, say ${\bf t}$ and ${\bf f}$.
Then, as discussed, for example, in \cite{M81} and  \cite{BM90}, letting $\rhd_{pc}$ be the restriction of $\rhd_c$ obtained when only the application of  propositional rules is  allowed, the rules for $\et$, the introduction rules for $\vel$ and $\imp$, together with the relations $A,\neg A\rhd_{pc}B$, $A, \neg B \rhd_{pc}\neg (A\imp B)$, determine the classical truth table for $\et$, half of the truth table for $\neg$ and three-fourth of the truth tables for $\vel $ and $\imp$. On the ground of the further relations $\neg A,\neg B \rhd_{pc} \neg (A\vel B)$ and $\neg A \rhd_{pc} A \imp B$, it then suffices to assume that $\neg A$ takes  the value ${\bf t}$, whenever $A$ takes the value ${\bf f}$, to obtain the classical truth tables.\footnote{Notice that none of the rules and relations concerning $\rhd_{pc}$ which are being used 
is specific to classical logic.}
When it comes to quantifiers we have that $v(F\{x/a\})$ ($F\{x/a\}$ pure), for $a$ an individual parameter,  has to take the value $\bf t$, whenever $v(\forall x F)$ takes the value $\bf t$, because of the $\forall$-elimination rule. Similarly $v(\exists x F)$ has to take the value $\bf t$, if for some parameter $a$, $v(F\{x/a\})$ takes the value $\bf t$, because of the $\exists$-introduction rule. As we will  show, an appropriate solution to our problem is obtained by simply reversing the last two implications, namely stating  that it is sufficient, for $v(\forall  x F)$ to take the value $\bf t$, that for every individual parameter $a$  of the language, $v(F\{x/a\})$ takes the value $\bf t$.  And similarly that it is necessary for $v(\exists x F)$ to take the value $\bf t$, that for some parameter $a$, $v(F\{x/a\})$ takes the value $\bf t$.

\subsection{Truth-value valuations}

\begin{definition}

 Let  ${\cal L}$ be a first order language. 
A truth-value  valuation ($tv$-valuation for short) of  ${\cal L}$ is a total  function $v$  from the collection
of pure atomic  formulae of ${\cal L}$ into $\{\bf{t}, {\bf f}\}$ such that $v(\bot) ={\bf f }$.

\end{definition}

A $tv$-valuation $v$ of ${\cal L}$ determines a unique extension $\bar{ v}$ to the pure formulae of ${\cal L}$, 
according to the classical  two-valued truth tables and the conditions:

\bi
\item ${\bar v}(\forall x H) ={\bf t}$ if and only it for every parameter $a$,  ${\bar v}(H\{x/a\}) ={\bf t}$.
\item  ${\bar v}(\exists  x H) ={\bf t}$ if and only it for some parameter $a$,  ${\bar v}(H\{x/a\}) ={\bf t}$.
\ei

\begin{definition}
$v$ $tv$-satisfies a pure formula $F$, if $\bar v(F)={\bf t}$; $F$ is $tv$-valid if every tv-valuation $v$ of ${\cal L}(F)$ satisfies $F$ and 
$F$ is a $tv$-semantic consequence of the pure formulae $G_1, \ldots, G_n, F$ if every $tv$-valuation $v$ of  ${\cal L}(G_1, \ldots, G_n, F)$, which $tv$-satisfies $G_1, \ldots, G_n$,  $tv$-satisfies $F$ as well.
\end{definition}

\subsection{Correctness and completeness for $tv$-semantics}

Correctness  and completeness   of the  $tv$-semantics  determined as  above by the $tv$-valuations,  for  the pure system $N_c$, holds.

\begin{theorem}
For $G_1, \ldots, G_n, F$ pure formulae, $G_1, \ldots, G_n \rhd_c F$ if  and only $F$ is a $tv$-semantic consequence of $G_1, \ldots, G_n$
\end{theorem}

{\bf Proof} Correctness is proved by a straightforward induction on the height of deductions in  pure $N_c$.
The only non entirely  trivial case occurs when the deduction ends with a $\forall:I$ or  $\exists:E$.
For example in the former case, letting  ${\cal D}$ be the immediate subderivation
with conclusion $H\{x/a\}$, given  any parameter $b$	 of ${\cal L}(G_1, \ldots, G_n, F)$, 
if $b$ is used as proper in (some  $\forall:I$ or $\exists:E$ rule applied in) $ {\cal D}$, 
we  first rename the occurrences of $b$ in ${\cal D}$ by a parameter $c$ new to ${\cal D}$
 and then replace  $a$ by $b$ throughout. 
 The result is a deduction of $H\{x/b\}$. By  the induction hypothesis, 
 any $tv$-valuation, which satisfies $G_1, \ldots, G_n$,  satisfies
 $H\{x/b\}$ as well. But that  means that it  satisfies $\forall x H$,
  as desired. Completeness can be proved, for example,  by applying the semantic tableaux method
   to pure formulae and considering only  parameters in the $\gamma$-reductions. If  $F$ is a consequence of $G_1, \ldots, G_n$,  the systematic tableaux procedure, initialized with $t.G_1, \ldots, t.G_n, f.F$, returns a closed tableaux from which a deduction ${\cal D}$ of $F$ from $G_1, \ldots, G_n$ can be obtained. Furthermore the  variables which have bound occurrences in ${\cal D}$ are exactly those which occur bound  in $G_1, \ldots, G_n, F$. $\Box$.
   
   \

 {\bf Note}  To have a correct and complete semantics  for  general formulae it suffices to state  that $F$ is a $tv$- semantic consequence
  of $G_1, \ldots, G_n$ if for some substitution
   $\theta=\{ x_1/a_1, \ldots, x_n/ a_n\}$, where
    $x_1, \ldots, x_n$ are the variables which have free occurrences in 
    $G_1, \ldots, G_n, F$, and $a_1, \ldots, a_n$
     are distinct parameters not occurring in 
     $G_1, \ldots, G_n,F$, we have that $F\theta$ is a pure semantic consequence  of
      $G_1\theta, \ldots, G_n\theta$. Correctness holds since from a deduction
       ${\cal D}$ of $F$ from $G_1, \ldots, G_n$, 
       after renaming the parameters among $a_1, \ldots, a_n$, 
       which  are used as proper in ${\cal D}$,
        one obtains a deduction of $F\theta$ from $G_1\theta, \ldots, G_n\theta$, 
        simply by replacing $x_1, \ldots, x_n$ by 
        $a_1, \ldots, a_n$ throughout ${\cal D}$. 
        As for completeness, we  first note that  its assumption and conclusion are invariant under renaming of bound variables. Therefore we may assume that no variable occurs both free and bound in $G_1, \ldots, G_n, F$. 
        Since, by assumption, $F\theta$ is a pure semantic consequence of $G_1\theta, \ldots, G_n\theta$, we may obtain a deduction of $F\theta$ from $G_1\theta, \ldots, G_n\theta$ in pure $N_c$,  which  is transformed into a deduction of $F$ from $G_1, \ldots, G_n$ simply by replacing $a_1, \ldots, a_n$ with $x_1, \ldots, x_n$ throughout.
  An immediate consequence is that the definition of $tv$-semantic consequence for general formulae does not depend on the choice of $\theta$.

\subsection{Equality}

Following $\cite{LH59}$, as  axioms for equality we take  reflexivity, namely $\forall(t=t)$, where $t$ is assumed to be parameter free and $\forall$ denotes  universal closure, and 
the  axiom of substitutivity of the form
\[
\forall(r=s \imp (F\{v/r\} \imp F\{v/s\}))
\]with $r, s$ and $F$ parameter free.
The two schemata of reflexivity and substitutivity will be denoted by $Rfl^{= s}$ 
and $Sbst^{= s}$.
$Rfl^{=s} $and 
$Sbst^{= s}$ 
are  easily seen to be equivalent over $N_c $
to $Rfl^{=s}$ and:
 
 \[
 \ba{ll}
Symm^{=s}& \forall (r=s \imp s=r),\\
Trans^{=s}& \forall(r=s\imp (s=t \imp r=t)),\\
 Cng^{= s}_p&\forall(r_1=s_1\et\ldots\et r_n=s_n \imp (p(r_1, \ldots, r_n) \imp p(s_1,\ldots, s_n)))\\
 Cng^{= s}_f& \forall(r_1=s_1\et\ldots\et r_n=s_n \imp 
  f(r_1, \ldots, r_n) = f(s_1,\ldots, s_n))
  \ea
  \]
   for any $n$-ary relation and function symbol $p$ and $f$, where all the terms shown are parameter free.
$N_c^=$ results from $N_c$ by allowing any formula in $Rfl^s$ and  $Sbst^s$ to be considered as a discharged assumption.

\

{\bf Note}  That  the equality axioms, formulated for variables only, namely $\forall x (x=x)$ and $\forall x \forall y (x=y \imp (F\{v/x\} \imp F\{v/y\})$,
are not sufficient for  a satisfactory development of the  logic of partial terms,
was first noticed in \cite{LH59}.

\subsection{$tv$-semantics for $N_c^=$}

\begin{definition}
A $tv$-valuation with equality of ${\cal L}$ is a $tv$-valuation of ${\cal L}$, which satisfies the axioms in $Rfl^{=s}$, $Symm^{=s}$, $Trans^{=s}$ and $Cng^{=s}$.
\end{definition}
 
 In other words,  $v$ is a $tv$-valuation  with equality  if  the binary relation
 $\{(r,s): v(r=s)={\bf t}\}$, to be denoted by $=^v$, 
 is a congruence relation with respect to the canonical interpretation
 of the function symbols $\{((t_1, \ldots, t_n), f(t_1, \ldots, t_n))\} $ and the relations $p^v =\{(t_1, \ldots, t_n): v(p(t_1, \ldots, t_n))={\bf t} \}$,  for $p$ relation symbol in ${\cal L}$, where $t_1, \ldots, t_n$ range over $PureTerm_{\cal L}$.

\ 

Correctness and completeness for $N_c^=$ holds with  respect to the notion of $tv$-semantic consequence based on $tv$-valuations with equality.

\begin{theorem}
For $G_1, \ldots, G_n, F$ pure formulae, $G_1, \ldots, G_n \rhd_c^= F$ if and only if  
every $tv$-valuation with equality of  ${\cal L}(G_1,\ldots, G_n, F)$ which $tv$-satisfies $G_1, \ldots, G_n$,  $tv$-satisfies $F$ as well.
\end{theorem}

{\bf Proof} Correctness is an immediate consequence of the correctness of $N_c$. Completeness can be achieved through the tableaux method, by interleaving 
the logical reduction steps with steps in which one appends, one after the other, the countably many judgments of the form $t.E$ where $E$ belongs to $Rfl^{=s}$, $Symm^{=s}$, $Trans^{=s}$ or $Cong^s$. $\Box$

\ 

Extension to general formulae can be obtained as for $N_c$.

\subsection{The extension property}

The following  property will be our  basic tool  for  dealing with $tv$-semantics for $N_c$ and $N_c^=$.

\begin{proposition}  Extension Property

If $v$ is a $tv$-valuation of ${\cal L}$ (with equality) and ${\cal L}\subset {\cal L}'$, then there is a map $\Phi$
from $PureTerm_{{\cal L}'}$ onto $PureTerm_{\cal L}$ and a valuation (with equality) $v'$ of ${\cal L}'$ such that:

\bi
\item[1)] for a term $t$ of ${\cal L} $ with variables among $x_1, \ldots, x_k$ and pure terms $r_1', \ldots, r_k'$ of ${\cal L}'$
\[
\Phi(t\{x_1/r_1', \ldots, x_k/r_k'\})= t \{ x_1/\Phi(r_1'), \ldots, x_k/\Phi(r_k') \},
\]
in particular if $t$ is a pure term of ${\cal L}$, $\Phi(t)=t$,
\item[2)] for  a  formula $F$ of ${\cal L}$ with free variables among $x_1, \ldots, x_k$ and pure terms $r_1', \ldots, r_k'$ of ${\cal L}'$
\[
\bar{v'}(F\{x_1/r_1', \ldots, x_k/r_k'\})= \bar{v} (F\{ x_1/\Phi(r_1'), \ldots, x_k/\Phi(r_k') \}),
\]
in particular if $F$ is a pure formula of ${\cal L}$,  then $\bar{v'}(F) = \bar{v}(F)$.

\ei
\end{proposition}

{\bf Proof}  For every $n$-ary function symbol  $f \in {\cal L}'\setminus {\cal L}$, fix a total function 

${\bf f}:PureTerm_{\cal L}^{~~n} \imp PureTerm_{\cal L}$ (for $n=0$, ${f}
$ is either  a constant or a parameter and ${\bf f}$ is a pure term, say $f_0$, of ${\cal L}$), which, in  case $v$ is a $tv$-valuation with equality,  is congruent with respect to $=^v$ (for example ${\bf f}$ can be any constant function). 
If $t$ is a parameter or a constant of ${\cal L}$,
 let $\Phi(t)=t$. If $t'$ is a parameter or a constant in ${\cal L}'\setminus {\cal L}$,
  let $\Phi(t')= t'_0$. 
  If $t'$ is $g(t_1', \ldots, t_n')$ with $g$
   in ${\cal L}$, let 
$\Phi(t')= g(\Phi(t_1') , \ldots, \Phi(t_n'))$, 
finally, if $t'$ is $f(t_1', \ldots, t_n')$ let 
$\Phi (t')= {\bf f} (\Phi(t_1'), \ldots, \Phi(t_n'))$. Furthermore for $p$,  $n$-ary relation symbol of ${\cal L}$ let
\[
v'(p(t_1', \ldots, t_n'))= v(p(\Phi(t_1'), \ldots, \Phi(t_n'))
\]
and,  for $q$ $n$-ary relation symbol in ${\cal L}'\setminus {\cal L}$,
let $v'(q(t_1', \ldots, t_n'))$ be defined arbitrarily provided $v'(q(s_1', \ldots, s_n'))={\bf t}$, whenever
 $v'(q(t_1', \ldots, t_n'))={\bf t}$ and 
 $v'(t_1'=s_1')={\bf t}, \ldots, v'(t_n'=s_n')={\bf t}$.
   
$1)$ and $2)$ are easily proved by induction on the height of $t$ and $F$ respectively. 
 $\Box$

\ 

Since, in the previous proof, it is the choice of ${\bf  f}$ which determines $v'$, we will say that $v'$ is the extension of $v$ based on ${\bf f}$.

\

{\bf Remark} The notion of $tv$-valuation can be relativized to any fixed subset ${\cal P}_0$ of the set of parameters of ${\cal L}$, 
assumed to be non empty, in case ${\cal L}$ has no constant,
by taking into account only the formulae whose parameters belong to ${\cal P}_0$ and considering only parameters in ${\cal P}_0$ in defining  the meaning of the quantifiers.  If ${\cal P}_0$ is infinite,  the proof of  correctness  remains unchanged. If ${\cal P}_0$ is finite,  correctness can be established along the lines of the previous proof. In fact  if $v_0$ is a
valuation restricted to any set of parameters ${\cal P}_0$ which satisfies $F$, then it suffices to note that $v_0$ can be extended to a valuation $v_0'$ of ${\cal L}(F)$, which still satisfies $F$,  by mapping   all the parameters which do non belong to ${\cal P}_0$ into any one of  the parameters in ${\cal P}_0$. 

Thus, for example, the $tv$-valuation restricted to 
$\{a, b\}$,

\[
\ba{rl}
v_0= & \{(p(a,a),{\bf t}), (p(b,b), {\bf t}), (p(a,b),{\bf f}), (p(b,a), {\bf f}), ((a=a),{\bf t}), \\
&((b=b),{\bf t}), ((a=b),{\bf f}), ((b=a),{\bf f})\}
\ea
\]
 which satisfies $\forall x \exists y p(x,y)$, 
but does not satisfy $\exists x \forall y p(x,y)$, 
suffices to show that in $N_c^=$ 
one cannot deduce the latter sentence from  the former. 
Similarly the $tv$-valuation restricted to $\{a\}$,

\[
\{ (p(c), {\bf t}), (p(a), {\bf f}), (a=a, {\bf t}), (c=c, {\bf t}), (a=c, {\bf f }), (c=a, {\bf f})\},\]
 for $c$ a constant, suffices to show that in $N_c^=$,  $\exists x p(x)$ cannot be deduced  from $p(c)$, and  the $tv$-valuation restricted to $\{a\}$:

\[
\ba{l}
 \{ (p(f^n(a), f^{n+1}(a)), {\bf t}): n\in N\} \cup \{(p(f^n(a),f^m(a)), {\bf f}):  m\neq n+1\}\\
 \cup  \{(f^n(a)=f^n(a) , {\bf t}): n\in N\} \cup \{(f^n(a)=f^m(a), {\bf f}): n \neq m \},
 \ea
 \]
  where $f^0(a)$ denotes $a$ itself,
suffices to show that  $\forall x\exists y p(x,y)$ is not deducible from $\forall x p(x,f(x))$.
On the other hand completeness for $tv$-valuations restricted to finite sets of parameters clearly fails.
For example $\exists x p(x,x)$ is not derivable in $N_c^=$ from $\forall x\exists y p(x,y)$ and $\forall x\forall y \forall z (p(x,y) \et p(y,z)\imp p(x.z)$,
 although it is satisfied by any $tv$-valuation   restricted to a finite set of parameters, which satisfies the latter two sentences.

\subsection{Totally denoting valuations}

{\bf Notation}  $t\downarrow$ denotes the formula $\exists y ~y=t$, for $y$ any variable not occurring in $t$.

The usual natural deduction system with equality, in which $\forall$-elimination and $\exists$-introduction
can be applied to any substitutable term, is easily seen to be equivalent to $N_c^=$, provided
$\forall(t\downarrow)$  is  allowed as a discharged assumption, for any term $t$. 
We denote with  $N_c^{\downarrow =}$ the resulting deduction system.  $N_c^{\downarrow =}$ is clearly equivalent to  $N^{=}_c$,   provided  formulae of the form $c\downarrow$  and $\forall x_1, \ldots, x_n f(x_1,\ldots, x_n)\downarrow$, for  all the constant $c$ and function symbol $f$ of the language, are allowed as discharged assumptions.

\begin{definition}
A tv-valuation with equality $v$ for ${\cal L}$ is said to be totally denoting if for every  pure term $t$ of ${\cal L}$,
$v$ tv-satisfies $t \downarrow $, 
namely there is a parameter
 $a$ such that $v(a=t) ={\bf t}$.
\end{definition}

\begin{proposition}
A tv-valuation $v$ for ${\cal L}$ with equality  is totally denoting if and only if every constant of ${\cal L}$ is denoting,
and for every $n$-ary function symbol $f$ and $n$-tuple of parameters $a_1, \ldots, a_n$, $f(a_1,\ldots, a_n)$ is denoting.
\end{proposition}

{\bf Proof} By a straightforward induction on the height of terms. $\Box$
\begin{theorem}
Correctness and completeness for $N_c^{\downarrow =}$ holds with  respect to the notion of tv-semantic consequence based on  totally denoting tv-valuations.
\end{theorem}

{\bf Proof} Immediate from the above propositions $\Box$.

\ 

{\bf Note} 
To every totally denoting $tv$-valuation $v$ for ${\cal L}$ there corresponds an elementarily equivalent set  theoretic interpretation
$I_v$.
The domain $D^{I_v}$ of $I_v$ is the set of parameters of ${\cal L}$. The interpretation of a constant  symbols in $I_v$ is 
a parameter $a$, such that $v(a=c)={\bf t}$. Similarly the interpretation of an $n$-ary function symbol $f$ is
a total function:
\[
f^{I_v}= \{((a_1, \ldots. a_n), b): v(b=f(a_1, \ldots, a_n))={\bf t}\}.
\]  Finally, for any relation symbol $p$ of ${\cal L}$, 
\[
p^{I_v} =\{ (a_1, \ldots, a_n): v(p(a_1, \ldots, a_n))={\bf t} \}.
\]
 Let $\tau$ be any assignment of elements of $D^{I_v}$ to variables and parameters which leaves all the parameters fixed, so that,  under $\tau$,  the value of any pure term $t$ is $t$ itself.
A  straightforward induction shows that if $F$ is a pure formula of ${\cal L}$, then $\bar{v}(F)={\bf t}$ if and only if $I_v, \tau \models F$. As a consequence for every sentence $F$ of ${\cal L}$,  $\bar{v}(F)={\bf t}$ if and only if $I_v \models F$, which is what we mean by saying that  $v$ and $I_v$ are elementarily equivalent. The  quotient of $I_v$ with respect to $=^v$ is a normal structure elementarily equivalent to $I_v$, therefore to $v$.  The  completeness theorem for (the ordinary set theoretic semantics) of $N_c^{\downarrow =}$ is thus an immediate consequence of the completeness of $tv$-semantics with equality
for $N_c^{\downarrow =}$.

\ 
\begin{proposition}
The Extension Property holds also for the totally denoting valuations.
\end{proposition}

{\bf Proof} If $v$ is totally denoting and $v'$
is an extension of $v$ to ${\cal L}'$, 
then $v'$ is also totally denoting since 
$\bar{ v'}(\exists x (x=t'))= \bar{v}(\exists x (x=\Phi(t'))$ and $\bar{v}(\exists x (x=\Phi(t')) = {\bf t}$,
 because $\Phi(t')$ is a pure term of ${\cal L}$ and $v$ is totally denoting. $\Box$

\subsection{Introducing the  undefined  $\uparrow$}

\begin{proposition}
 A totally denoting tv-valuation $v$  of ${\cal L}$ can  be extended to a tv-valuation $v^{\uparrow}$ with equality of the language ${\cal L} + \uparrow$, where $\uparrow$ is a constant not belonging to ${\cal L}$, such that for every pure formula $F$ of ${\cal L}$, $\bar{v}(F) = \bar{v^\uparrow} (F)$ and $\uparrow$ is non denoting with respect to $v^\uparrow$.
 \end{proposition}
  
{\bf Proof}  We set $v^\uparrow(r=s)= {\bf t}$ if and only if $r=s$ belongs to the 
 smallest set of equalities between pure terms of  ${\cal L} +\uparrow$, which contains all the equalities $t'=t'$ and $r=s$ such that $v(r=s)= {\bf t}$ and furthermore
 contains $f(r_1, \ldots, r_n) = f(s_1, \ldots, s_n)$ whenever  for all $1\leq i\leq n$ it already contains $r_i=s_i$.
On  all the remaining pure atomic formulae which contain $\uparrow$, $v^\uparrow$ takes the value
${\bf f }$ and  $v^\uparrow(A)= v(A)$ for every pure  atomic formula of  ${\cal L}$. The claim follows by a straightforward induction on the height of $F$. To prove that $v^\uparrow$ is a valuation with equality it suffices to show that $v^\uparrow(r_1=s_1\et \ldots \et r_n=s_n \imp (p(r_1, \ldots, r_n)\imp p(s_1, \ldots, s_n)))={\bf t}$. If all of $r_1,\ldots, s_n$ belong to ${\cal L}$, that holds since $v^\uparrow$ agrees with $v$, which is a $tv$-valuation with equality. Thus let us assume that, for example $\uparrow$ occurs in $s_i$. Then, by definition,  $v^\uparrow(p(s_1, \ldots, s_n))={\bf f}$, and we have to show that also $v^\uparrow(p(r_1, \ldots, r_n))={\bf f}$. That follows from the fact that if $\uparrow$ occurs in $s_i$ and $v^\uparrow(r_i=s_i)={\bf t}$, then $\uparrow$ occurs also in $r_i$. As a matter of fact we have that 
if $v^\uparrow(r=s)={\bf t}$ and $\uparrow$ occurs in $s$ then $\uparrow$ occurs also in $r$ and conversely, as it  follows immediately from the definition of $v^\uparrow$ on equalities.
Obviously  that guaranties also that  $\uparrow$ cannot be denoting. $\Box$

  \subsection{Strictness}
  \begin{definition}
  Let $N_c^{=s}$ be the result of adding to $N_c^=$ the following {\em strictness}  axioms:
  
  \bi
  \item[1)] $c\downarrow$
  \item[2)] $\forall(f(t_1, \ldots, t_n)\downarrow \imp t_1\downarrow\et\ldots\et t_n \downarrow)$
  \item[3)] $\forall(p(t_1, \ldots, t_n)\imp t_1\downarrow \et \ldots \et t_n\downarrow)$ for every relation symbol $p$ other than $=$, and 
$t_1,\ldots, t_n$  parameter free.
   \ei
  A strict tv-valuation  of ${\cal L}$ is a tv-valuation of ${\cal L}$ with equality  which satisfies the strictness
  axioms.
  \end{definition}
  
  In $3)$ we have to leave aside $=$, since otherwise, from the adoption of $t=t$ as an axiom, it would follow that every $t$ is defined. Thus  our notion of strictness is more relaxed  than the one  usually adopted when the existence predicate is taken as primitive (see \cite{F95} for example).
  
  \ 
  
   The proof of  correctness and completeness of the semantics based on totally denoting $tv$-valuations for $N_c^=$ can be easily adapted to establish the following:

  \begin{theorem}
Correctness and completeness for $N_c^{= s}$ holds with  respect to the notion of tv-semantic consequence based on  strict  tv-valuations.
\end{theorem}
 
\begin{proposition}
The Extension Property holds also for  strict  tv-valuations, provided the extension  is based on functions ${\bf  f}$ which are strict, namely satisfy the following condition:

 $a)$ if $v({\bf f}(r_1, \ldots, r_n)\downarrow)={\bf t}$, then $v(r_1\downarrow)={\bf t}, \ldots, v(r_n \downarrow)={\bf t}$.
\end{proposition}

{\bf Proof} If $v$ is strict and $v'$ is an extension of $v$ to ${\cal L}'$ based on
a function ${\bf f}$ satisfying condition $a)$,
 then $v'$ is also  strict.
For, assume ${\bar v'} (f(t_1', \ldots, t_n')\downarrow)={\bf t}$, namely  ${\bar v'} (\exists x (x=f(t_1', \ldots, t_n'))={\bf t}$. If $f\in {\cal L}'\setminus {\cal L}$ by the Extension Property it follows
that ${\bar v}(\exists x (x= {\bf f }(\Phi(t_1'), \ldots, \Phi(t_n')))={\bf t}$. 
By the strictness of ${\bf f}$, it follows that $\Phi(t_1')\downarrow, \ldots, \Phi(t_n')\downarrow$, 
 namely ${\bar v}(\exists x_1(x_1=\Phi(t_1')))={\bf t}, \ldots, {\bar v}(\exists x_n(x_1=\Phi(t_n')))={\bf t}$, from which, by the Extension Property again, we may conclude that
  $\bar{v'}(\exists x_1(x_1=t_1'))={\bf t}, \ldots, \bar{v'}(\exists x_n (x_n =t_n'))={\bf t}$, 
  namely $\bar{v'}(t_1'\downarrow)={\bf t}, \ldots, \bar{v'}(t_n'\downarrow)={\bf t}$, 
  as required for $v'$ to be strict. The case in which $f\in {\cal L}$ or 
${\bar v'}(p(t_1', \ldots, t_n'))={\bf t}$, for $p$ other than $=$, is entirely similar. $\Box$

\

 {\bf Note}
As for totally denoting $tv$-valuations, to every strict valuation $v$ of ${\cal L}$ there corresponds an elementarily equivalent  (partial) set theoretic interpretation $I_v$ of ${\cal L}$. $D^{I_v}$ is still the set of  parameters of ${\cal L}$ but $f^{I_v}$ is, in general,  a partial function.  For a given  assignment $\sigma$ of elements of $D^{I_v}$ to variables and parameters,
the value $\sigma(t)$ which  $t$ takes under $\sigma$ is an element of $D^{I_v}$ iff  $t\sigma$ is a denoting term, namely $v(t\sigma\downarrow)={\bf t}$.  $I,\sigma\models F$ is defined by letting $I, \sigma\models r=s$ iff $v(\sigma(r),\sigma(s))={\bf t}$ (even if $\sigma(r)$ or $\sigma(s)$ does not belong to $D^{I_v}$); for $p$ other that $=$, $I_v, \sigma \models p(t_1, \ldots, t_n)$ if and only if                                     
 $\sigma(t_1), \ldots, \sigma(t_n)$ belong to $D^{I_v}$ and $(\sigma(t_1), \ldots, \sigma(t_n)) \in p^{I_v}$
 (namely $v(p(\sigma(t_1), \ldots, \sigma(t_n))= {\bf t}$)). For compound formulae $I_v, \sigma \models F$ is defined as usual.  For every pure formula $F$  of ${\cal L}$ and assignment $\tau$, which leaves  the parameters fixed, ${\bar v}(F)={\bf t}$ if and only if $I_v, \tau \models F$, so that for a sentence $F$,  ${\bar v}(F)={\bf t}$ if and only if $I_v \models F$. 
For,  if  $F$ is of the form $p(t_1, \ldots t_n)$, from ${\bar v}(F)={\bf t}$, by the strictness of $v$, it follows that $t_1, \ldots, t_n$ are all denoting terms, so that $\tau(t_1), \ldots, \tau(t_n)$ belong to $D^{I_v}$, and $(\tau(t_1), \ldots, \tau(t_n)) \in p^{I_v}$ so that $I_v, \tau \models F$.
 As a consequence we have the completeness of $N_c^{=s}$ with respect to partial set theoretic interpretations.

\ 
 
 {\bf Note}  If  a $tv$-valuation $v$   is extended into $v^\downarrow$, rather than into  $\bar{v}$, by using the clauses:
\bi
\item[a)] ${v^\downarrow }(\forall x H) ={\bf t}$ iff for every pure term $t$,  ${v^\downarrow}(H\{x/t\}) ={\bf t}$.
\item[b)]  ${ v^\downarrow }(\exists  x H) ={\bf t}$ iff for some pure term $t$,  ${v^\downarrow}(H\{x/t\}) ={\bf t}$.
\ei
then  a straightforward modification of the previous arguments shows that the resulting semantics is correct and complete with respect to the usual natural deduction system, without equality, in which $\forall$-elimination and $\exists $-introduction can be applied to any substitutable term, do be denoted by 
$N_c^\downarrow$, and that  the Extension Property  still holds.
Furthermore $v^\downarrow$ is elementarily equivalent to a (total) set theoretic structure $I_{v^\downarrow}$, whose  domain is the set $D^{I_v}$ of the pure terms of the language, so that the usual completeness theorem  for $N_c^\downarrow$ immediately follows.
The same applies if $v$ is  a valuation with equality, thus obtaining a correct and complete semantics
for $N_c^{\downarrow =}$. Since if $v$ is a totally denoting valuation, then obviously $\bar v=v^\downarrow$, the $tv$- semantics for $N_c^\downarrow$ based on $v^\downarrow$ subsumes the one based on totally denoting $tv$-valuations, so that its completeness can also be inferred from the completeness of the latter.  As in the previous case, one can also immediately infer the usual completeness theorem for $N_c^{\downarrow =}$.
That shows   the interest of $tv$-semantics  even  if one is concerned only with  {\em total} classical logic with or without equality.
In particular  the standard classical set theoretic semantics can be rather  effectively introduced as a very natural generalization of $tv$-totally denoting semantics, by replacing the fixed domain of the pure terms of the language  by an arbitrary non empty set and the {\em total} canonical interpretation of the function symbols by their interpretation with arbitrary total functions on such a set.
Concerning  the last point,  we wish to note the difficulty  one  faces in motivating
the choice of totality, if the classical set theoretic structures are to be presented as a model of ordinary mathematical structures, which may carry partial, rather than only total, operations, like the reals with the $x^{-1}$ or  $log$ function, for example.

\section{Conservativeness of partial selection functions}

\begin{theorem}
If $D$ is  a formula of ${\cal L}$ with distinct free variables $x_1, \ldots, x_n,y$, and $f$ is an  $n$-ary function symbol not in ${\cal L}$, then the conjunction of the following two formulae is conservative over ${\cal L}$ with respect to $N_c^=$:
\[
\ba{l}
\epsilon^1_y(f; D)~~~~\forall(f(x_1, \ldots, x_n)\downarrow \imp \exists y D)\\
\\
\epsilon^2_y(f;D)~~~~\forall(\exists y D \imp \exists y (y=f(x_1,\ldots, x_n)\et D))
\ea
\]namely if $G_1, \ldots, G_n, F$ are formulae of ${\cal L}$, $f$ does not occur in $G_1,\ldots, G_n, F$ and $G_1, \ldots, G_n, \epsilon^1_y(f; D),\epsilon^2_y(f;D) \rhd_c^= F$, then  $G_1, \ldots, G_n \rhd_c^= F$.
The same holds for $N_c^{=s}$.
\end{theorem}

{\bf Proof}  We deal first with the case in which $G_1, \ldots, G_n, F$ are pure.  By the correctness and completeness of the $tv$-semantics with equality for $N_c^=$,  it suffices to show that  the Extension Property can be applied to any
 $tv$-valuation with equality $v$ of ${\cal L}$,   so as to obtain a valuation $v'$ of ${\cal L}+f$ which  satisfies $\epsilon^1_y(f; D)$ and  $\epsilon^2_y(f; D)$.
  If $v$ is not totally denoting, fix a non denoting term $t_0$ of ${\cal L}$
 and an enumeration of all the parameters of ${\cal L}$. 
 If $t_1,\ldots, t_n$ are all denoting terms of ${\cal L}$,  $a_i$ is the first parameter in the fixed enumeration such that $v(a_i=t_i)={\bf t}$ and $b$ the first one such that
 $v(D\{x_1/a_1,\ldots, x_n/a_n, y/b\})={\bf t}$, provided there is such a $b$, we let
 ${\bf f}(t_1, \ldots, t_n)= b$; if on the contrary there is no $b$ such  that  $v(D\{x_1/a_1,\ldots, x_n/a_n, y/b\})={\bf t}$ or for some $1\leq i\leq n$, $t_i$ is non denoting, then we let ${\bf f}(t_1, \ldots, t_n)=t_0$.
 As it is easy to check,  ${\bf f}$ is congruent with respect to $=^v$, so that the extension $v'$ of $v$ to ${\cal L}+f$, based on ${\bf f}$,  is a $tv$-valuation with equality, and it  is also strict.
  Furthermore $\overline{v'}$ satisfies $\epsilon^1_y(f;D)$ and $\epsilon^2_y(f;D)$.
  Since $\epsilon^1_y(f;D)$ follows in $N^=_c$ from $\forall x_1\ldots\forall x_n\forall y (f(x_1,\ldots, x_n)=y \imp D)$, it suffices to verify that $\overline{v'}$ satisfies the last formula, namely that, for every $n+1$-tuple of parameters $a_1,\ldots, a_n, b$, if 
  $  \overline{v'}(f(a_1, \ldots, a_n)=b$, then $  D\{x_1/a_1,\ldots, x_n/a_n,y/b\})={\bf t}
  $
 By  the  Extension Property $\overline{v'}(f(a_1, \ldots, a_n)=b)=\overline{v}({\bf f }(a_1, \ldots, a_n) =b)$. Thus from  $\overline{v'}(f(a_1, \ldots, a_n)=b)={\bf t}$ it follows that $\overline{v}({\bf f }(a_1, \ldots, a_n)=b)={\bf t}$, which, by the definition of ${\bf f}$, it can only happen if $\overline{v}(D\{x_1/a_1,\ldots, x_n/a_n,y/b\})={\bf t}$.
  
 As for $\epsilon^2_y(f;D)$, we have  to verify that for every $n$-tuple of parameters $a_1,\ldots, a_n$,  
 if $\overline{v'}(\exists y D\{x_1/a_1,\ldots,x_n/a_n\})={\bf t}$, then 
  $\overline{v'}(\exists y(y=f(a_1, \ldots,a_n)\et D\{x_1/a_1,\ldots, x_n/a_n\}))={\bf t}$.
   From the assumption, by the Extension Property   it follows that
  $\overline{v}(\exists y D\{x_1/a_1,\ldots,x_n/a_n\})={\bf t}$. Thus there is a parameter $b$, which we may assume is the first in the given enumeration, such that 
  $\overline{v}(D\{x_1/a_1,\ldots,x_n/a_n, y/b \})={\bf t}$. Therefore ${\bf f }(a_1,\ldots, a_n)=b$.
On the other hand,  $\overline{v'}(\exists y(y=f(a_1, \ldots,a_n)\et D\{x_1/a_1,\ldots, x_n/a_n\}))={\bf t}$ if and only if 
    there is a parameter $c$ such that 
$ \overline{v'}(c=f(a_1, \ldots,a_n) \et D\{x_1/a_1,\ldots, x_n/a_n, y/c\})={\bf t}.
   $
  By the Extension Property that holds if and only if 
   there is a parameter $c$ such that
  $\overline{v}(c ={\bf f}(a_1, \ldots, a_n) \et D\{x_1/a_1,\ldots, x_n/a_n, y/c\})={\bf t}
  $.
 Therefore  it suffices to take $b$ for $c$  to conclude that our claim holds.
   If  $v$ is totally denoting, it suffices to consider its extension with the "undefinite" $v^\uparrow$
  and replace $t_0$ by $\uparrow$ in the previous argument,  to obtain the desired extension  of $v$.
   By the Extension Property for strict valuation the result  applies to $N_c^{=s}$ as well.
  To extend the result to general formulae it suffices to repeat the argument given  for the extension of the completeness theorem.  $\Box$

\subsection{Conservativity of  partial  description  functions}

\begin{theorem}
If $D$ is  a formula of ${\cal L}$ with distinct free variables $x_1, \ldots, x_n,y$, and $f$ is an  $n$-ary function symbol not in ${\cal L}$, then the following  formula is conservative over ${\cal L}$ with respect to $N_c^=$:

\[
\iota_y(f;D)~~~ \forall( f(x_1, \ldots, x_n)=y \equiv D\et \forall y' (D\{y/y'\} \imp y'=y))
\]
The same holds for $N_c^{=s}$.i
\end{theorem}

{\bf Proof} Given $D$, let $D^!$ be $D\et \forall y' (D\{y/y'\} \imp y'=y)$. By the proof of the first part of the previous theorem applied to $D^!$, we can conservatively add the formula 
$a)~~\forall (f(x_1,\ldots, x_n)=y \imp D^!)$.
Furthermore we can conservatively add $\epsilon^2_y(f;D^!)$.
From $D^!$ it logically follows $\exists y(D\et \forall y'(D\{y/y'\} \imp y'=y))$, from which by $\epsilon^2_y(f;D^!)$ it follows $\exists y (y=f(x_1,\ldots, x_n)\et D\et \forall y'(D\{y/y'\} \imp y'=y))$. Let then $z$ be such that $z=f(x_1,\ldots, x_n)\et D\{y/z\}\et \forall y'(D\{y/y'\} \imp y'=z)$.  From $z=f(x_1,\ldots, x_n)\et D\{y/z\}$ by $D^!$ it follows that $z=y$ hence $f(x_1,\ldots, x_n)=y$. Thus also the reverse implication in $a)$, and therefore $\iota_y(f;D)$, is deducible in the conservative extension provided by the previous theorem with respect to $D^!$. Hence  $\iota_y(f;D)$ is conservative over ${\cal L}$ with respect to $N_c^=$. $\Box$

\begin{corollary}
Under the assumption of the previous theorem 
\[
\forall(f(x_1, \ldots, x_n) =y \equiv D)
\]
is conservative over  $N_c^= + U_y D$ where $U_y D$ states the uniqueness condition for $y$ satisfying $D$, namely $\forall (D \et D\{y/y' \} \imp y'=y)$.
\end{corollary}

{\bf Proof} Under $U_y D$,  $D\et \forall y' (D\{y/y'\} \equiv y'=y)$ and $D$ are obviously  logically equivalent, so that it suffices to substitute the latter for the former in
$\iota_y(f: D) $, in the previous theorem. $\Box$

\

{\bf Directions for further work} 
As we noticed, the notion of strictness we have adopted is tailored to fit the proposal in \cite{LH59}, to deal with singular terms, hence 
in  doesn't  assume that if $t=t$ holds, then $t$ is denoting. It would be interesting to match the present treatment with the  more demanding notion of strictness, by finding an appropriate axiomatization of equality. The $tv$-semantic approach to the conservativity of partial description functions and of  partial selection functions, in the latter case under the assumption of the  {\em determinacy} of equality, namely the assumption  $\forall x\forall y (x=y \vel x\neq y)$, should be extended to the case of intuitionistic logic.
Obviously such questions call also for a proof theoretic treatment. That requires a preliminary investigation of logic with equality  and the proof of an appropriate  subterm and subformula property (for cut free derivations in a suitable  sequent calculus). Joint work with F. Previale in that  direction  is well under way.

\subsection{Acknowledgment}

We wish to express our thanks to  Alberto Marcone and Flavio Previale for helpful conversations and remarks.

  \end{document}